\documentclass[12pt]{amsart}
\newtheorem{theorem}{Theorem}[section]
\newtheorem{lemma}[theorem]{Lemma}
\newtheorem{proposition}[theorem]{Proposition}

\newtheorem{example}[theorem]{Example}
\newtheorem{definition}[theorem]{Definition}
\newtheorem{corollary}[theorem]{Corollary}

\newcommand{\N}{\mbox{$\mathbb{N}$}}

\newcommand{\K}{\mbox{$\mathbb{K}$}}

\newcommand{\E}{\mbox{${\mathcal E}$}}

\newcommand{\R}{\mbox{${\mathcal R}$}}

\newcommand{\Li}{\mbox{${\mathcal L}$}}
\newcommand{\M}{\mbox{${\mathcal M}$}}

\newcommand{\T}{\mbox{${\mathcal T}$}}

\begin{document}
   \title[One-sided ideals and approximate identities in operator algebras]
{One-sided ideals and approximate identities in operator algebras}
\author{David P. Blecher*} 

\thanks{*This research was 
supported by a grant from the National Science Foundation}  
\thanks{*Many of these results were presented 
at the November '01 Meeting of the A.M.S. in Irvine, California.}
\thanks{Sept 14, 2001.  Revision of November 12, '01.}  
\address{Department of Mathematics\\University of Houston\\
4800 Calhoun\\Houston, TX 77204-3008, U.S.A.} 
\email{dblecher@math.uh.edu}\maketitle

\vspace{10 mm}
 
\begin{abstract}
A left ideal of any $C^*$-algebra is an example of an 
operator algebra with a right contractive approximate 
identity (r.c.a.i.).   
Indeed 
left ideals in $C^*$-algebras may be characterized 
as the class of such operator algebras, which happen 
also to be triple systems.
 Conversely, we show here and in 
a sequel to this paper \cite{BK}, that 
operator algebras with r.c.a.i. should be studied in
terms of a certain left ideal of a $C^*$-algebra.   
We study left ideals from the perspective of `Hamana
theory' and using the multiplier algebras introduced by 
the author.  
More generally, we
develop some general theory for operator algebras
which have a 1-sided identity or approximate identity,
including a Banach-Stone theorem for these algebras,
and an analysis of the `multiplier operator algebra'.
\end{abstract}

\pagebreak
\newpage

\section{Introduction and notation}
A left ideal of any $C^*$-algebra is an example of an
operator algebra with a right contractive approximate
identity.   
Conversely, we shall see here, and in a 
sequel to this paper \cite{BK}, that
operator algebras with 
a right contractive approximate  
identity  should be studied in
terms of a certain left ideal of a $C^*$-algebra. 
 
A (concrete) operator algebra is a subalgebra of $B(H)$,
for some Hilbert space $H$.   More abstractly,
an operator algebra will be an algebra $A$ with
a norm defined on the space $M_n(A)$ of $n \times n$ 
matrices with entries in $A$, for each $n \in \N$, such that
there exists a {\em completely isometric}\footnote{A map 
$T : X \rightarrow Y$ is 
completely isometric if $[x_{ij}] \mapsto
[T(x_{ij})]$ is isometric on $M_n(X)$ for all $n \in \N$.} 
homomorphism
$A \rightarrow B(H)$ for some Hilbert space $H$.
In this paper all our operator algebras and spaces will be taken
to be complete.  We shall say that an 
operator algebra is {\em unital} if it has a two-sided
contractive identity.   Unital operator algebras were 
characterized abstractly in \cite{BRS}.  However the class of 
nonselfadjoint operator algebras which is of most interest
to $C^*$-algebraists or those interested in noncommutative 
geometry is the class of left or right ideals in a 
$C^*$-algebra, which
is a nonselfadjoint operator algebra with a 
right (respectively left) contractive approximate identity.   
Thus in this paper we shall study operator algebras
with a one-sided (usually right) contractive approximate
identity. 
We shall 
abbreviate `right (resp. left) contractive approximate  identity'
to `r.c.a.i.' (resp. `l.c.a.i.').  It is well known even for 
Banach algebras that an algebra with both a left and a right
c.a.i. has a two-sided c.a.i..
  
As we point out in \S 2, 
left ideals in $C^*$-algebras may be characterized
as the class of nonselfadjoint 
operator algebras with r.c.a.i, which happen
also to be `triple systems'.   This therefore suggests
that left ideals in $C^*$-algebras may profitably 
be studied using machinery that exploits both the
`operator algebra' and the `triple' structure, and indeed 
we do take this approach here.   For example `morphisms'
of left ideals will be what we call `ideal homomorphisms'
below, namely homomorphisms which are also 
`triple morphisms'. 
 
Much of this paper is concerned with left ideals
in unital operator algebras which have
a r.c.a.i..    In a formal sense this coincides with the 
class of all operator algebras with r.c.a.i., since every
operator algebra is easily seen to be a left ideal in some
operator algebra.  By Proposition 6.4 in
\cite{BEZ} (although we shall not use this fact except for
motivational purposes), this also coincides with the 
class of left M-ideals in
unital operator algebras.   We shall not however use 
this last fact, except as evidence
that these objects are more important than they seem at first.
We will make the blanket convention in 
this paper
that all ideals, left or otherwise, are closed, i.e. complete.

We are not aware of any general results in the literature on 
operator algebras with r.c.a.i.,
and one of the main goals of this paper and 
its sequel \cite{BK} is to take on such a study.
Although operator algebras with r.c.a.i. have not hitherto come up
in the literature on general 
operator algebras as much as their two-sided relative,
 they do arise very naturally
(e.g. left ideals in a $C^*$-algebra,
$Ae$ for a projection $e$ in a
unital operator algebra $A$;
or Example \ref{rowf}, or
Theorem 4.11 in \cite{BMP}).   
We also felt that it was 
worthwhile to dispel the common feeling that only 
operator algebras with 2-sided approximate identities have 
any satisfactory theory.   Indeed we show amongst other
things that they 
have an abstract characterization, Banach-Stone type theorems,
reasonable multiplier algebras (which are 
operator algebras with two sided identity of 
norm 1), and they have an
operator space predual if and only if they are
dual operator algebras' in the usual strong sense 
of that term.   Also, this subject becomes a little 
more interesting with a certain `transference principle' 
in mind.  This principle (which is proved fully in 
the sequel \cite{BK}),
allows one to deduce many general results about 
 operator algebras  with r.c.a.i., from results about
left ideals in a $C^*$-algebra.  Namely there is an important 
left ideal ${\mathfrak J}_e(A)$ of a $C^*$-algebra $\E(A)$, which 
is associated to any 
such operator algebra.  We call ${\mathfrak J}_e(A)$ the `left ideal 
envelope' of $A$.   This is 
 analoguous to what 
happens in the case of operator algebras with 2-sided c.a.i., 
which are largely studied these days
in terms of a certain $C^*$-algebra, namely
the $C^*-$envelope.   

 \vspace{4 mm}

We now describe the layout of the paper.  
In \S 2 we give many preliminary lemmas which will be used  
later and in \cite{BK}, as well as many new results about left ideals and
operator algebras with r.c.a.i..

In \S 3 we look at Banach-Stone type theorems.        
The classical Banach-Stone theorem (see e.g. \cite{Con} IV.2)
may be stated in 
the following form: if $C(K_1) \cong C(K_2)$ linearly 
isometrically, then they are *-isomorphic (from 
which it is clear that the compact spaces
$K_1$ and $K_2$ are homeomorphic).  
Indeed the usual proofs show that the linear  
isometry equals a *-isomorphism $C(K_1) \rightarrow C(K_2)$
multiplied by a fixed unitary in $C(K_2)$.   Note that
the linear isometry is unital (i.e. takes the $1 \mapsto 1$)
if and only if $u = 1$, and one often proves this special 
case first.
There are numerous noncommutative versions of this, the most 
well known due to Kadison \cite{Ka}, where 
the $C(K)$ spaces are replaced by $C^*$-algebras.  For unital 
maps between nonselfadjoint operator algebras with 
two-sided identities of norm 1, a similar theorem is true,
that is,  any unital linear complete isometry between
such operator algebras, is also a
homomorphism.  This 
 is an immediate consequence of the existence of the
Arveson-Hamana $C^*-$envelope of \cite{Arv1,Arv2,Ham1}
(indeed this generalization of
the Banach-Stone theorem was a major consideration in
\cite{Arv1,Arv2}, see also \cite{Bcomm,ERns}).
In  \cite{BShi} Appendix B.1 it is  
shown that a linear
surjective complete isometry $\varphi : A \rightarrow B$ between
operator algebras with contractive approximate identities
may be written as $\varphi = \pi(\cdot) u$ for a
surjective completely isometric homomorphism 
$\pi : A \rightarrow B$, and a unitary $u$ with $u, u^* \in M(B)$.
Here $M(B)$ is the multiplier algebra (see \cite{PuR} for example). 
In \S 3 below we  examine 
such theorems in the case of operator algebras with 
one-sided identities or contractive
approximate identities, for example in the case of
left ideals in a 
$C^*$-algebra.   We see that the Murray-von Neumann
equivalence of projections plays a role here.   

In \S 4 we study the `left multiplier operator algebra' $LM(A)$
of an  operator algebra $A$ with l.c.a.i. (which will be 
a symmetrical theory to that of $RM(A)$ for an 
operator algebra $A$ with r.c.a.i).     
In the next paper we develop a satisfactory candidate for  
$LM(A)$ if $A$ has a r.c.a.i. - this does not work out 
quite as nicely as the case considered
in \S 4, unless $A$ is a
left ideal of a $C^*$-algebra.  

We remark that we are not aware even
of much general theory of one-sided
contractive approximate identities in general
Banach algebras.  The only references we can give for 
this subject are 
the works of P. G. Dixon (see \cite{Pal}
for references), and the general texts \cite{BonsallDuncan,Pal}.

Often it is convenient to state only the `right-handed' version of a
result, say.  For example Lemma \ref{namo} is a result about
operator algebras
with r.c.a.i..   Of course by symmetry
there will be a matching `left-handed' version, in our example
it will be about operator algebras with l.c.a.i..  If we want to
invoke this `left-handed' version, we will refer to the
`other-handed  version of Lemma \ref{namo}', for example.  Sometimes we will be
sloppy and simply refer to Lemma \ref{namo}, and leave it to the reader
to check that the result referred to does indeed have a symmetric
or `other-handed  version'. 

We end the introduction with some
more notation, and some background results which will
be useful in various places.  
We reserve the letters $H, K$ for Hilbert spaces,
and $J$ for a left ideal of a $C^*$-algebra.

If $A$ is an algebra then we write $\lambda
: A \rightarrow Lin(A)$ for the  canonical `left regular
representation' of $A$ on itself.
If $S$ is a 
subalgebra of $A$ then 
the left idealizer of $S$ is the 
subalgebra $\{ x \in A : x S \subset S \}$ of $A$.
Note $S$ is a left ideal in this subalgebra, whence the name.
Similarly for the right idealizer; the (2-sided) 
idealizer is the intersection of the left and right idealizer.

By a `representation' $\pi : A \rightarrow B(H)$ of 
an operator algebra $A$ we shall
mean a completely contractive homomorphism.  If
 $A$ has r.c.a.i. and if we say that 
$\pi$ is {\em nondegenerate}, then at the very least  
we mean that 
$\pi(A) H$ is dense in $H$.  Note that this last 
condition does not imply
in general
 that $\pi(e_\alpha) \zeta \rightarrow \zeta$ for $\zeta \in H$,
where $\{ e_\alpha \}$ is the r.c.a.i.,
as one is used to in the two-sided case.  One also cannot appeal to
Cohen's factorization theorem in its usual form
(see however \cite{Pal} \S 5.2).    

We will use without comment several very basic facts from
$C^*$-algebra theory (see e.g. \cite{Ped}), such as
the basic definitions
of the left multiplier algebra $LM(A)$ of a $C^*$-algebra,
and the multiplier algebra $M(A)$.

As a general reference for operator spaces the reader might consult
\cite{ERbook,Pis} or \cite{P}.
We write $\; \hat{ } \; : X \rightarrow X^{**}$ for the 
canonical map, this is a complete isometry if $X$ is an
operator space, and is a homomorphism if $X$ is an operator 
algebra (giving the second 
dual the Arens product \cite{BonsallDuncan}).

We will often consider the basic examples $C_n$ (resp. $R_n$)
of operator algebras with right (resp. left) identity of norm 
1; namely the $n \times n$ matrices
`supported on' the first column (resp. row).  This is a
left (resp. right) ideal of $M_n$, and has the projection
$E_{11}$ as the 1-sided identity.  
We write $C_n(X)$ for 
the first column on $M_n(X)$, that is
$M_{n,1}(X)$.  If $X$ is an operator space
so is $C_n(X)$.    

If $X$ and $Y$ are subsets of an operator algebra we usually write 
$X Y$ for the {\em norm closure} of the set of finite 
sums of products $x y$ of a term in $X$ and a term in $Y$.
For example, if $J$ is a left ideal of a $C^*$-algebra $A$,
then with this convention $J^* J$ will be
 a norm closed $C^*$-algebra.  This convention extends to 
three sets, thus $J J^* J = J$ for a left ideal of a $C^*$-algebra
as is well known (or use the proof of Lemma \ref{lidci} below 
to see this).   
 
We recall more generally   that a TRO ({\em ternary ring of 
operators}) is a (norm closed for this paper) subspace $X$
of $B(K,H)$ such that $X X^* X \subset X$.   It is well known
(copy the proof of Lemma \ref{lidci} below) that in this case 
$X X^* X = X$.  Then $X X^*$ and $X^* X$ are $C^*$-algebras,
which we will call the left and right $C^*$-algebras of $X$ respectively,
and $X$ is a $(X X^*) - (X^* X)$-bimodule.
It is also  well known that TRO's are 
the same thing as Hilbert $C^*$-modules.  A linear map
$T : X \rightarrow Y$ between TRO's is a {\em triple 
morphism} if $T(x y^* z) = T(x) T(y)^* T(z)$ for all 
$x,y,z \in X$.   TRO's are operator spaces, and 
triple morphisms are completely contractive, and indeed are 
completely isometric if they are 1-1 (see e.g. \cite{Ham2},
this is related to results of Harris and Kaup).
A completely
isometric surjection between TRO's is a triple morphism.  This 
last result might date back to around 1986, to Hamana and 
Ruan's PhD thesis independently.
See \cite{Ham2} or \cite{BShi} A.5 for a proof.

If $J$ is a left ideal in a $C^*$-algebra, and if $B$ is a
$C^*$-algebra, then
a certain class of maps $J \rightarrow B$ are 
of particular importance in this paper.
Namely, those which
are restrictions of
*-representations $\theta : J J^* \rightarrow B(H)$.  We will
call these {\em ideal homomorphisms} of $J$.  They are 
characterized in \ref{babi} below as the 
homomorphisms which are also triple morphisms.

Next we recall 
the left multiplier algebra
$\M_\ell(X)$ of an operator space $X$.   This is a unital
operator algebra, which may be viewed as a subalgebra of 
$CB(X)$ containing $Id_X$, but with a different (bigger in general)
norm.   
Here $CB(X)$ are the `completely bounded' linear maps
on $X$.   
One may take the definition of $\M_\ell(X)$ from the following
result from \cite{BEZ}:
 
\begin{theorem} \label{bez}  A linear 
$T : X \rightarrow X$ on an operator space is in 
$Ball(\M_\ell(X))$ if and only if $T \oplus Id$ is a
complete contraction $C_2(X) \rightarrow C_2(X)$.
\end{theorem}

As we said above, $\M_\ell(X)$ is an operator algebra.  The
matrix norms on $\M_\ell(X)$ may be described via the 
natural isomorphism $M_n(\M_\ell(X)) \cong \M_\ell(M_n(X))$.
That is the norm of a matrix $[T_{ij}]$ of multipliers may be taken to 
be the norm in $\M_\ell(M_n(X))$ of the map 
$$[x_{ij}] \mapsto [\sum_k T_{ik}(x_{kj})] \; \; .$$

There are several other
equivalent definitions of $\M_\ell(X)$ given in
\cite{BShi,BEZ,BPnew}.  One of our main motivations 
for the introduction of multipliers of 
 operator spaces in \cite{BShi}
was in order to give 
a more unified and `extremal' approach to the 
main theorems characterizing operator algebras and 
modules.   We pointed out in \S 5 of \cite{BShi} 
that in order to prove the characterization of 
operator algebras of \cite{BEZ} say, it is clearly only 
necessary to check that the `left regular representation'
$\lambda : A \rightarrow CB(A)$, is a complete 
isometry into the operator algebra
$\M_\ell(A)$.  But this is immediate
from a theorem such as \ref{bez} above
- see the simple proof
of Lemma \ref{lemu} for more details if required.
(This deduction of \cite{BEZ} from \ref{bez} was
noticed independently by Paulsen, who also gave a simplified 
proof of \ref{bez}).   The same proof works to
give immediately
the following characterization of operator algebras with
one-sided c.a.i., which is a slight variation on
 \cite{BMDO} 1.11 to which the reader is referred if 
more details of proof are needed (we warn the
reader that the statement of that result
is stated slightly incorrectly, it neglected to mention
that the matrices in $A$ there should have norm $\leq 1$).

\begin{theorem}  \label{1sBRS}  Let $A$ be an operator space
which is an algebra with a right identity of norm 1
or r.c.a.i..  Then $A$
is completely isometrically isomorphic to a concrete operator
algebra (via a homomorphism of course), 
if and only if we have
$$ \Vert (x \oplus Id_n) y \Vert \leq 1$$
for all $n \in \N$ and
$x \in Ball(M_n(A)), y \in Ball(M_{2n,n}(A))$.
\end{theorem}
 
To explain the notation of the theorem, we have written
$Id$ for a formal identity,
thus the expression $(x \oplus Id_n) y$ 
above means that the
upper $n \times n$-submatrix of $y$ is left
multiplied by $x$, and the lower submatrix is left alone.

{\em Acknowledgements:}  We thank M. Kaneda for catching 
several misprints, and for several conversations
related to the present paper and \cite{BK}.  
 
\section{General results and lemmata}

In this section we collect many simple
lemmas and other background facts which will be used later
and in \cite{BK}.  Some of the results here for 
left ideals in $C^*$-algebras may be known, but in any case
our proofs are short.
Certainly the first lemma is classical and 
extremely well known; we have appended a (well known)
proof which might be convenient for the reader who accepts
the existence of c.a.i.'s in $C^*$-algebra 
(modern texts
in this area contain short proofs of this
latter fact (see e.g. \cite{Fill})).

\begin{lemma} \label{lidci} (Classical) 
A (norm closed) left ideal $J$ in a $C^*$-algebra is an 
operator algebra with a positive right contractive approximate
identity.   Also $J \cap J^* = J^* J \subset J \subset J J^*$,
 so that $J$ is a left ideal of $J J^*$. 
\end{lemma}

\begin{proof}  A left ideal $J$ in a $C^*$-algebra $A$ is
clearly a subalgebra of $A$.  Also $J J^*$ and $J^* J$ are
 $C^*-$subalgebras of $A$.  So $J^* J$ has a positive 
c.a.i. $\{ e_\alpha \}$.  Then for $x \in J$,  
$$\Vert x e_\alpha - x \Vert^2 =
\Vert e_\alpha x^* x e_\alpha - x^* x e_\alpha - 
e_\alpha x^* x + x^* x \Vert \rightarrow 0 \; \; .$$
The remaining assertions follow  immediately
from this; for example if $x \in J \cap J^*$ then
$x^* = \lim x^* e_\alpha$, so that $x \in J^* J$.
\end{proof}

\vspace{3 mm}

Note that $J^* J$ also equals 
$\{ x \in J : e_\alpha x \rightarrow x \}$, where 
$\{ e_\alpha \}$ is the c.a.i. for $J$ mentioned above.  This is part 
of our motivation for the next definition and result (which will
only be used in part of \S 4)\footnote{Another part of our motivation 
and inspiration is the theory in
 \cite{Bhmo} and \S 4 of \cite{BMP}.}.

\begin{definition} \label{defp}  We say that an operator algebra $A$
 with r.c.a.i. (resp. l.c.a.i.) has property ($\R$) (resp. ($\Li$)) 
if a  r.c.a.i. (resp. l.c.a.i.) $\{ e_\alpha \}$ 
exists for $A$ such that $ e_\alpha e_{\alpha'} \rightarrow 
e_{\alpha'}$ 
(resp. $e_{\alpha'} e_\alpha \rightarrow e_{\alpha'}$)
for each fixed $e_{\alpha'}$ in the net.    In 
this case we define $\R(A) = \{ x \in A : e_\alpha x \rightarrow 
x \}$ (resp. $\Li(A) = \{ x \in A : x e_\alpha \rightarrow x \}$).
    \end{definition}

\noindent 
{\bf Remarks.}  We note that a left ideal 
of a $C^*$-algebra has property 
($\R$), and in this case $\R(A) = J^* J$.   More generally
a subalgebra of a $C^*$-algebra with a self-adjoint  
right c.a.i. has property ($\R$), since in this case 
$(e_\alpha e_{\alpha'})^* = e_{\alpha'} e_\alpha
\rightarrow e_{\alpha'} = e_{\alpha'}^*$.
An operator algebra
with two-sided c.a.i. obviously
has property ($\R$), and in this case $\R(A) =
A$.   Certainly every operator algebra with a
right identity of norm 1 has property ($\R$).
We are not aware
of any  operator algebras  with r.c.a.i. which do not have 
property ($\R$) (if all do this would certainly solve a
problem encountered in \cite{BMP} \S 4 and elsewhere
in some of the authors work: whether every `Hilbertian
module' is a rigged module).  Thus it seems a 
reasonable and general property.    

\begin{proposition} \label{defpt}   If an operator algebra $A$ 
with r.c.a.i. has  property ($\R$), 
then $\R(A)$  is a norm closed right ideal of 
$A$ (and hence is an operator algebra)
 with two sided c.a.i..  Moreover, $\R(A)$
does not depend on the particular 
c.a.i. $\{  e_\alpha \}$ considered.   Also,
$A \; \R(A) = A$ and $\R(A) \; A = \R(A)$.
 
Similar results hold 
for property ($\Li$).
\end{proposition}

\begin{proof}  The first assertion we leave as a simple exercise.
Suppose that $A$ has property ($\R$) with respect to one
r.c.a.i. $\{  e_\alpha \}$, and let $\{ f_\beta \} $ be another
r.c.a.i. such that $f_\beta f_{\beta'} \rightarrow  f_{\beta'}$ 
for every fixed $\beta'$.   Let $B = \{ a \in A :
f_\beta a \rightarrow a \}$, another right ideal of $A$ with 
two sided c.a.i..  Note that $\R(A) B = \R(A)$ and $B \R(A) = B$.
Thus by (the other-handed version of)
 \cite{BMP} Theorem 4.15, $B = \R(A)$.
The remaining assertions are left to the reader.
\end{proof}

\begin{example} \label{rowf} \end{example}
Let $B$ be a unital operator algebra, a unital subalgebra
of a $W^*$-algebra $N$, 
and define $M_\infty(N)$ to be the von Neumann algebra
$B(\ell^2) \bar{\otimes} N$, thought of as infinite matrices 
$[b_{ij}]$  with entries $b_{ij}$ indexed over $i, j \in \N$.
We let $M_\infty(B)$ be the subset of $M_\infty(N)$ consisting 
of those matrices with entries $b_{ij}$ in $B$.
Let $C^w_{\infty}(B)$ be the `first column' of $M_\infty(B)$,
and let $R_\infty(B)$ be the space of row vectors 
$[ b_1 b_2 \cdots ]$ with  entries  $b_i \in B$, such that 
$\sum_k b_k b_k^*$ converges in norm.   These spaces are 
familiar objects in operator space theory (see \cite{ERbook}), and 
also in $C^*$-module theory.
We may then consider the closed subspace 
$A = C^w_{\infty}(B) R_\infty(B)$ of
$M_\infty(B)$;  those familiar with operator space theory will have 
no trouble verifying that $A$ is a subalgebra of 
$M_\infty(N)$, that $A$ has a nonnegative r.c.a.i., and 
indeed if $B = N$ then $A$ is
a left ideal of $M_\infty(B)$.
   In fact  
$A$ contains the $C^*$-algebra $\K_\infty(B)$, namely 
the spatial tensor product $\K(\ell_2) \otimes B$ (which in 
the language of $C^*$-modules equals $\K(C_\infty(A))$), 
and the usual c.a.i. for this $C^*$-algebra, namely 
$I_n \otimes 1_B$, is a r.c.a.i. for $A$.  Thus $A$ has property
($\R$).  It is easily verified 
that $\K_\infty(B)$ is a right ideal in $A$, and in fact 
$\R(A) = \K_\infty(B)$.   This algebra $A$ is important to those 
working on the borders of operator spaces and operator algebras,
although it has rarely appeared in the literature as far as 
we are aware, and that fleetingly.   

\vspace{4 mm}

The next lemma concerns `principal ideals'.
By a `principal ideal' in a $C^*$-algebra $A$, we will mean 
by analogy with pure algebra,
 an ideal of the form $Ax$  for some $x \in A$.
We are not taking the norm closure here,
$Ax = \{ a x : a \in A \}$ for some $x \in A$; however 
in view of the
importance of closed ideals in $C^*$-algebra theory, below we
will only consider principal ideals which are already norm
closed. 
We remark however that if $A$ is 
nonunital, then when transferring results of 
the kind found in pure ring theory
to the $C^*$-algebra case, one finds that the 
multiplier algebra of $A$ plays 
an important role.  Hence the reader might argue that
ideals of the form $Ax$ with $x \in M(A)$ (or $RM(A)$), 
should be thought of as  `principal ideals' of $A$ too.   
However in  our 
paper we shall not need these more general 
ideals, and so stick with the earlier definition 
for convenience.
Interestingly though,
all these kind of ideals have a simplified form:

\begin{proposition} \label{prid}  Let $A$ be a $C^*$-algebra, 
and suppose that $J = Ax$ is a closed left ideal, with $x \in A$ (resp. 
$x \in M(A)$).  Then $J = Ae$ where $e$ is an orthogonal 
projection in 
$J$ (resp. in $M(A)$).  
\end{proposition}

\begin{proof}  Since $J$ is the range of an adjointable map on $A$,
$J$ is orthogonally complemented in the sense of 
$C^*-$module theory, by \cite{W-O} 15.3.9.
This implies that $J = Ae$ where $e$ is an orthogonal
projection in $M(A)$.  So that if $A$ is unital we are done,
and note that in this case $Ae$ has a right identity of 
norm $1$.
However in any case, if $x \in A$, then $A x = M(A) x$
 (clearly $A x \subset M(A) x$, but if $T \in M(A)$ 
then $T x = \lim T e_\alpha x \in A x$).   Thus
applying the above we see that $J$ has a right identity $f$
of norm $1$, and $f \in J \subset A$.  Hence $J = Af$.
\end{proof}

\begin{corollary} \label{babe}
Suppose that $J$ is a left ideal of a $C^*$-algebra, and suppose
that $J$ has a right identity.  Then 
$J$ also has a right identity of norm 1.  Moreover
the latter is the norm limit of any r.c.a.i. for the left ideal.
\end{corollary}

\begin{proof}   Let $J$ be the left ideal, which by the
previous result, has a right identity $e$ of norm 1.
So $e = e^* \in J \cap J^* = J^* J$.  If $\{ e_\alpha \}$ is a r.c.a.i. for 
$J$ then $\{ e_\alpha^* e_\alpha \}$ is a 2-sided c.a.i. for 
$J^* J$ (see e.g. \ref{pid2}), thus $
e_\alpha^* e_\alpha = e_\alpha^* e_\alpha e \rightarrow e$.  Finally,
$$\Vert e_\alpha - e \Vert^2 = \Vert e_\alpha^* e_\alpha
- e_\alpha^* e - e e_\alpha + e \Vert \rightarrow 0 \; \; . $$
\end{proof} 

\vspace{4 mm}

Later we
 will prove the analogous result to the last corollary valid 
for operator algebras.  

\begin{lemma} \label{pid}  Suppose that $A$ is an operator
algebra with two right identities $e$ and $f$ of norm 1.
Then $e = f$.
\end{lemma}
 
\begin{proof}   Since $e$ and $f$ are orthogonal
projections we have $e = e f = e^* = f e = f$.
\end{proof}

\begin{lemma} \label{namo}  (\cite{BonsallDuncan} 28.7).
If $A$ is an operator algebra with r.c.a.i.
then $A^{**}$ is an operator algebra with right 
identity of norm 1.  If $A$ has a right identity $e$,
then $\hat{e}$ is the right identity of $A^{**}$.
\end{lemma}

\vspace{4 mm}
      
The last 
assertion of the following 
lemma is well known; however we state 
it here because it is a trivial corollary of the more 
important (for us) first part:

\begin{lemma} \label{pid2} Suppose that $a \in B(H)$, and 
$\{ e_\alpha \}$ is a net of contractions in $B(H) $ such that
$a e_\alpha \rightarrow a$.  Then
$a e_\alpha e_\alpha^* \rightarrow a$, 
$a  e_\alpha^* e_\alpha \rightarrow a$, and 
$a  e_\alpha^* \rightarrow a$.

Thus if $J$ is a left ideal of a $C^*$-algebra, and 
if $\{ e_\alpha \}$ is a r.c.a.i. for $J$,
  then $\{ e_\alpha^* e_\alpha \}$ is a nonnegative 
right contractive
approximate identity for $J$ (and indeed also is a
2-sided c.a.i. for 
the  $C^*$-subalgebra $J \cap J^* = J^* J$).
\end{lemma}

\begin{proof}   We use a technique from \cite{BMP}.
If $a e_\alpha \rightarrow a$
then $a e_\alpha e_\alpha^* a^* \rightarrow a a^*$, so 
that $0 \leq a (I - e_\alpha e_\alpha^*) a^*  \rightarrow 0$.
Thus by the $C^*$-identity, 
$a \sqrt{I - e_\alpha e_\alpha^*} \rightarrow 0$.  Multiplying 
by $\sqrt{I - e_\alpha e_\alpha^*}$ we see that 
$a (I - e_\alpha e_\alpha^*)  \rightarrow 0$ as required
for the first assertion.  Also, $$\Vert a e_\alpha^* - a \Vert 
\leq \Vert a e_\alpha^* - a e_\alpha e_\alpha^* \Vert +
\Vert a e_\alpha e_\alpha^* - a \Vert \;  \rightarrow 0$$
since  
$\Vert a e_\alpha^* - a e_\alpha e_\alpha^* \Vert  
\leq \Vert a - a e_\alpha \Vert \rightarrow 0$.  
Finally, $$\Vert a e_\alpha^* e_\alpha - a \Vert
\leq \Vert a e_\alpha^* e_\alpha - a e_\alpha \Vert 
+ \Vert a e_\alpha - a \Vert \leq
\Vert a e_\alpha^* - a \Vert +  \Vert a e_\alpha - a \Vert
\rightarrow 0$$
by what we just proved.
\end{proof}
 
\vspace{4 mm}

Note that the last result shows that a r.c.a.i. for a  $C^*$-algebra
is a l.c.a.i. too.

If $J$ is a left ideal in a $C^*$-algebra, then 
we recall from \S 1 that an {\em ideal representation} 
or {\em ideal homomorphism} of $J$ is a
restriction of                                                                
a *-representation $\theta : J J^* \rightarrow B(H)$ to $J$.
Clearly such a map is completely contractive.  

\begin{proposition}  \label{babi}  Let $J$ be a left ideal of a
$C^*$-algebra, and let $\pi :  J
\rightarrow B(H)$ be a function.  Then $\pi$ is the restriction
of a *-representation $\theta : J J^* \rightarrow B(H)$
if and only if $\pi$ is a homomorphism and a 
triple morphism.  Moreover such $\pi$ is completely isometric if and only if
$\pi$ is 1-1, and  if and only if $\theta$ is 1-1.
\end{proposition}

\begin{proof}   If $\pi$ is the restriction
of a *-representation then it is evident that 
$\pi$ is a homomorphism and a
triple morphism.   Conversely, it is well known (see \cite{Ham2}),
that if  $\pi$ is a triple morphism, then there is an associated 
*-homomorphism $\theta : J J^* \rightarrow B(H)$ with the 
property that $\theta(xy^*) = \pi(x) \pi(y)^*$ for all 
$x, y \in J$.  If in addition $\pi$ is a homomorphism, and 
$\{ e_\alpha \}$ is a 
positive r.c.a.i. for $J$, then $\{ \pi(e_\alpha) \}$
is a positive r.c.a.i. for $\pi(J)$, and so for $x \in J$,  
$$\theta(x) = \lim \theta(x e_\alpha) = \lim \pi(x) \pi(e_\alpha)^*
= \pi(x)$$
by Lemma \ref{pid2}.  

If further $\pi$ is 1-1, then it is shown in \cite{Ham2} that 
$\theta$ is 1-1.
\end{proof}

\vspace{ 4 mm}

The following result is a simple consequence of the fact that
$J J^* J = J$:
 
\begin{lemma} \label{trmo}  Let $J$ be a left ideal of a
$C^*$-algebra, and let $ \theta : J J^* \rightarrow B(H)$ be a
 *-homomorphism.  If $\pi$ is the restriction to $J$ then
$\theta $ is nondegenerate if and only if
 $\pi(J) H$ is dense in $H$.
\end{lemma}

If $A$ has a left identity of norm 1 but no right identity,
and if $\pi : A \rightarrow B(H)$ is a nondegenerate  
isometric representation, then $\pi(e) = Id$, 
so that $\pi(ae) = \pi(a)$, so that 
$ae = a$ for all $a \in A$.  This is a contradiction.
Thus there is in general little point in seeking nondegenerate
isometric representations of algebras with l.c.a.i..
(The reader may think at this juncture of `adjoint nondegeneracy'
but this really is a different issue to the point we
are now making).
   The following discussion, proposition, and 
subsequent definition of $\Li$-isometric representation,
gives one way to fix the above problem.

If $A$ has left identity $e$ of norm 1, then $A$ clearly has
property ($\Li$) of \ref{defp}, and this identity
is the 2-sided identity  of $\Li(A) = A e$.  Moreover, the map
$A \rightarrow \Li(A)$ taking $a \mapsto a e$, is a completely
contractive homomorphism, and also is a complete quotient map and
indeed is a projection onto $\Li(A)$.   On the other hand,
if  $A$ has a l.c.a.i. and property ($\Li$),
 then we can make similar assertions
for the second dual using \ref{namo}: there is a completely
contractive homomorphism, which is a complete quotient map and
indeed a projection $A^{**} \rightarrow \Li(A)^{**}$.  
This is the map $F \mapsto FE$, where $E$ is a weak* limit 
point of the c.a.i. of $\Li(A)$.    
 
\begin{proposition} \label{reps}
Suppose that $A$ is an operator algebra with l.c.a.i. and
property ($\Li$) of \ref{defp}. 
 Let $\pi : A \rightarrow B(H)$ be a completely contractive 
representation
(resp. and also $\pi(A) H$ is dense in $H$).
Then $\pi_{|_{\Li(A)}} : \Li(A)  \rightarrow B(H)$ is a
completely contractive homomorphism (resp. and also
such that $\pi(\Li(A)) H$ is dense in $H$).  Conversely, if
$\theta : \Li(A)  \rightarrow B(H)$ is a completely contractive homomorphism,
then there exists a completely contractive homomorphism
$\pi : A \rightarrow B(H)$ extending $\theta$.  If further
$\theta(A) H$ is dense in $H$ then $\pi$ is unique, and
$\pi(A) H$ is dense in $H$.  Finally,
$$\{ T \in B(H) : T \pi(A) \subset \pi(A) \} 
= \{ T \in B(H) : T \pi(\Li(A)) \subset \pi(\Li(A)) \} \; . $$
\end{proposition}
 
\begin{proof}   The first statements are simple
exercises.  For the converse,
given such $\theta  : \Li(A)  \rightarrow B(H)$, consider the
series of completely contractive homomorphisms
$$A \hookrightarrow A^{**} \rightarrow \Li(A)^{**} 
\overset{\theta^{**}}{\rightarrow} B(H)^{**}
\rightarrow B(H) . $$
The homomorphism $A^{**} \rightarrow \Li(A)^{**}$ is the map described
above the Proposition, and the other maps are the canonical ones.
The composition of these homomorphisms is the desired $\pi$.  We 
leave it to the reader to check the details.
  Since $\pi(a) \theta(b) \zeta = \pi(ab) \zeta = \theta(ab)  \zeta$
for $a \in A, b \in  \Li(A), \zeta \in H$ we see that $\pi$ is
unique if $\pi(A) H$ is dense.

Finally, using the `other-handed version' of the last assertion of
\ref{defpt}, we see for example that if $T \pi(A) \subset \pi(A)$
then $T \pi(\Li(A)) = T \pi(A) \pi(\Li(A)) \subset
\pi(A) \pi(\Li(A)) = \pi(\Li(A)).$   The other
direction is similar.   \end{proof}                   

\vspace{4 mm} 
      
The previous result shows that $A$ and $\Li(A)$ have the same 
representation theory.  Thus the following definition which 
plays a role in \S 4 is 
somewhat natural: we say that a nondegenerate representation
$\pi : A \rightarrow B(H)$ is completely `$\Li$-isometric', if 
$\pi_{|_{\Li(A)}}$ is completely isometric on $\Li(A)$.  

\begin{lemma}  \label{lemu}  Let $A$ be an operator algebra
with a  r.c.a.i..  Then the canonical `left regular 
representation' of $A$ on itself yields completely contractive
embeddings (i.e. 1-1 homomorphisms)  
$$A \hookrightarrow \M_\ell(A) \hookrightarrow CB(A) \; \; ,$$
and the first of these embeddings, and their composition,
are completely isometric.   
\end{lemma}

\begin{proof}   Let $\lambda : A \rightarrow CB(A)$ be the
left regular representation.   This map is certainly 
completely contractive, however since $\lambda(a)(e_\alpha) 
= a e_\alpha \rightarrow a$ it is clear that it is a 
complete isometry.  It is also clear that if 
$a \in Ball(A)$ then $\lambda(a)$ satisfies the criterion
of \ref{bez},
so that $\lambda(a) \in Ball(\M_\ell(A))$.
A similar argument works at the matrix level. 
Thus $\lambda$ factors through $\M_\ell(A)$ via the
two completely contractive homomorphisms above.  Since 
$\lambda$ is completely isometric, so is the first embedding.
\end{proof}

\vspace{4 mm}       

We remark that the canonical
 embedding $\M_\ell(A) \hookrightarrow CB(A)$,
where $A$ is an operator algebra with r.c.a.i.,
 is not in general completely isometric, or even isometric
(an example of this is given in \cite{BK}).   
This has implications for our theory 
of multipliers in \cite{BK}.

\vspace{4 mm}

The following
`Hamana-theory' type results, which are very interesting 
in their own right, will
be one of our main tools to deduce results 
about operator algebras with r.c.a.i., from results
about left ideals in a $C^*$-algebra.  We give condensed versions 
of two of these results from \cite{BK}.

We say that a pair $(J,i)$ consisting 
of a left ideal $J$ in a $C^*$-algebra, and a completely isometric
homomorphism $i : A \rightarrow J$, is 
a {\em left ideal extension} of $A$
if $i(A)$ `generates $J$ as a TRO'.  That is, the span of 
expressions of the form $i(a_1) i(a_2)^* i(a_3)  i(a_4)^* 
\cdots i(a_{2n+1})$,
for $a_i \in A$, are dense in $J$.   It follows from this
that $\{ i(e_\alpha) \}$ is a r.c.a.i. for $J$ if 
$\{ e_\alpha \}$ is a r.c.a.i. for $A$.  

\begin{theorem} \label{upie}  Let $A$ be an 
operator algebra  with r.c.a.i..  Then there 
exists a left ideal extension $({\mathfrak J}_e(A),j)$ of $A$, 
with ${\mathfrak J}_e(A)$ a left ideal in a $C^*-$algebra
$\E(A)$, such that for 
any other left ideal extension $(J,i)$ of $A$, there exists 
a (necessarily unique and surjective)
ideal homomorphism $\tau: J 
\rightarrow {\mathfrak J}_e(A)$
 such that $\tau \circ i = j$.  Thus 
${\mathfrak J}_e(A)/(ker \; \tau)  \cong J$ completely isometrically
homomorphically (i.e as operator algebras) too.
Moreover $({\mathfrak J}_e(A),j)$ is unique in the following sense:
given any other $(J',j')$ with this universal property,
then there exists a surjective completely isometric 
homomorphism $\theta : {\mathfrak J}_e(A) \rightarrow J'$ such that 
$\theta \circ j = j'$.

Finally, $({\mathfrak J}_e(A),j)$ is a {\em triple envelope} for 
$A$ in the sense of \cite{Ham2}.
\end{theorem}

\vspace{4 mm}

We call $({\mathfrak J}_e(A),j)$ the {\em left ideal envelope} of 
$A$, and often write $\E(A) = {\mathfrak J}_e(A) {\mathfrak J}_e(A)^*$.
The map $j$ will be called the {\em Shilov embedding
homomorphism}.  From the last assertion of the theorem,
and the first  definition of $\M_\ell(A)$ given in 
\S 4 of \cite{BShi}, we may identify 
$\M_\ell(A)$ with $\{ R \in LM(\E(A)) : R j(A)  \subset j(A) \}.$
If $\E(A)$ is represented nondegenerately as a 
$C^*$-subalgebra of $B(H)$, then we can also identify
$\M_\ell(A)$ with  $\{ R \in B(H) : R j(A)  \subset j(A) \}$, 
completely isometrically isomorphically.  We shall not however
use this last remark.

\begin{corollary} \label{notfi}  Let $A$ be an operator
algebra with r.c.a.i., and $\lambda$ the usual left 
regular representation of $A$.  Then for any $T \in \M_\ell(A)$, 
regarded as a map on $A$,
we have $T \lambda(A) \subset \lambda(A)$.  Thus 
elements of
$\M_\ell(A)$, considered as maps on $A$,
are right $A$-module maps.  That is, $\M_\ell(A) \subset CB_A(A)$ as
sets.   Also, $\M_r(A) \subset \; _ACB(A)$ as sets.
\end{corollary}

\begin{proof}  The first assertion follows from the penultimate 
line before the statement of the Corollary, together with the 
fact that $j$ is a homomorphism.  For if $a \in A$, then the 
map $b \mapsto ab$ on $A$, corresponds to the map 
$j(b)  \mapsto j(a) j(b)$ on $j(A)$.  Thus if the left 
multiplier $T$ corresponds to an $R \in LM(\E(A))$ with 
  $R j(a) = j(T(a))$
then $R j(a b) = j(T(a) b)$ for any $b \in A$, which amounts to
the first assertion, and also yields the second assertion
immediately.   The third is similar. 
\end{proof}  

\vspace{5mm}

\begin{theorem} \label{Thlemi}  (\cite{BK}) 
Let $A$ be an operator algebra
with a r.c.a.i..   Then there exists an injective unital $C^*$-algebra 
$B$ containing $\E(A)$ as a $C^*$-subalgebra, and 
an orthogonal projection $e \in B$, such that $(B e,j)$ is an
injective envelope for $A$, where $j$ is  the 
Shilov embedding homomorphism.   
Moreover if $A$ has a right identity $f$ then $j(f) = e$.
\end{theorem}
 
We define the {\em ideal injective envelope} of $A$ 
to be the ideal $Be$ in the last theorem. 

\begin{corollary} \label{liie}  Let $J$ be a left ideal 
in a $C^*$-algebra.   Then 
the canonical embedding of a left ideal in a $C^*$-algebra,
into its ideal injective envelope $B e$ 
is an ideal homomorphism. 
\end{corollary}

\begin{proof}   Using \ref{upie}, there exists a
surjective ideal homomorphism $\tau : J \rightarrow {\mathfrak J}_e(J)$,
 such that $\tau = j$.  
Clearly ${\mathfrak J}_e(J)$ is a subalgebra and sub-TRO of 
$\E(J)$, and consequently also of $B$ and $B e$.  
\end{proof} 

\vspace{4 mm}

As another corollary of \cite{BK}, one may give 
a characterization of left ideals in $C^*$-algebras. 
We will  say that an operator space $X$ is an
abstract triple system if it is linearly completely isometrically
isomorphic to a TRO $Z$.  Note that then one may pull back the
triple product on $Z$ to a triple product $\{ \cdot ,
\cdot , \cdot \}$ on $X$, and by
the TRO result of Harris and Ruan mentioned in the introduction,
 this triple product on $X$ is unique,
i.e. independent of $Z$.  That is, this triple product is completely
determined by the `operator space structure' or matrix norms
on $X$.   Indeed Neal and Russo have a striking recent characterization
of abstract triple systems in terms of these matrix norms
\cite{NR}.  Putting this together with the theorem below, and
 a characterization of operator algebras with
right contractive approximate identity
(r.c.a.i.) which we gave in Theorem 1.2 for example, gives a
`completely abstract' characterization of 
left ideals in $C^*$-algebras.
 
\vspace{4 mm}
 
\begin{theorem}  \label{chli}  Let $A$ be an abstract operator algebra
which is also an abstract triple system (we are assuming the
underlying matrix norms for both structures coincide).
The following are equivalent:
\begin{itemize}
\item [(i)]  $A$ has a r.c.a.i. for the algebra product.
\item [(ii)]  $A$ has a r.c.a.i. $\{ e_\alpha \}$ for
the algebra product such that $\{ a , e_\alpha , b \}
\rightarrow ab$ for all $a,b \in A$;
\item [(iii)]  There exists a left ideal $J$  in a $C^*$-algebra,
and a surjective complete isometry
$A \rightarrow J$ which is both a
homomorphism (i.e. multiplicative), and a
triple morphism.
\end{itemize}
\end{theorem} 
\begin{proof}   It is trivial that (ii) implies (i),
and that (iii) implies (ii) (since we
can choose a nonnegative r.c.a.i. in $J$).  In fact it is not too
hard to see that (iii) implies (ii) for {\em any}
r.c.a.i. for $A$, but we will not need this.
 
Finally, given (i), we observe that by hypothesis,
$A$ is a a triple envelope of $A$.  Let ${\mathcal E}$ be
the `left  $C^*$-algebra' of this triple system, and let
$\langle \cdot | \cdot \rangle :
A \times A \rightarrow {\mathcal E}$ be the
associated sesquilinear inner product.
By \cite{BK}
Theorem 2.3, the map $\psi(x) = \lim_\alpha
\; \langle x  | e_\alpha \rangle$
is a well defined completely isometric triple morphism
of $A$ onto a left ideal $J$ of ${\mathcal E}$,
and $\psi = \psi \circ Id_A$ is also a homomorphism.
\end{proof}
 
\vspace{1 mm}
 
\noindent {\bf Remarks.}  (1).   There are various rather
trivial `TRO' characterizations of left ideals in $C^*$-algebras.
For example, one can characterize left ideals in
$C^*$-algebras with right identities as the TRO's $Z$
which possess an
element $e$ of norm 1 such that  $xe = x$ for all $x \in Z$.   This is
equivalent to $x e^* = x$ for all $x \in  Z$.  A similar condition with a
net does the general case.   Indeed a similar idea gives
a rather trivial characterization of  left ideals in $C^*$-algebras
as the `abstract triple systems' which are also operator
algebras with a r.c.a.i.  if one
adds a algebraic  compatibility conditions such as (ii)
above, between the operator algebra
product and the triple product.   Our result seems a little deeper than
this.  The only compatibility condition we seem to need
between the operator algebra
 ond the triple product is that the induced matrix norms are the
same.

\vspace{1 mm}
 
(2).   A slight modification of this result also gives a
characterization of $C^*$-algebras, by replacing
`r.c.a.i.' by `c.a.i.'.    We are grateful to Bernie Russo
for pointing out a recent paper \cite{ECP}
which gives such a characterization, but without needing the
matrix norms.

\vspace{1 mm}
 
(3)  It would be interesting if, in the spirit of
\cite{NR}, one could give a purely linear
characterization of left ideals in $C^*$-algebras.
There is such a result in \cite{BEZ}, but it makes reference
to the containing $C^*$-algebra in the hypotheses. 
 
\vspace{4 mm}

Corollary \ref{notfi} allows us to generalize the main 
result of \cite{BMDO} (see also \cite{LM3}) to 
algebras with one-sided c.a.i.:
 
\begin{theorem} \label{BLM1s}  Let $A$ be an operator algebra with
r.c.a.i., which has a predual operator space.  Then
$A$ has a right identity $e$ of norm 1.  Also $A$ is a
`dual operator algebra', which means that the product
on $A$ is separately weak* continuous, and there exists
a completely isometric homomorphism, which is also
a homeomorphism with respect to the weak* topologies, of $A$
onto a $\sigma$-weakly (i.e. weak*-) closed subalgebra $B$
of some $B(H)$.
\end{theorem}
 
\begin{proof}   The first assertion appears in \cite{BMDO}
Theorem 2.5 (indeed for this part
we only need a predual Banach space).  From \cite{BMDO} Theorem 3.2,
$\M_\ell(A)$ is a dual operator algebra.
We saw in \ref{lemu} and
\ref{notfi} that  $\lambda : A \rightarrow \M_\ell(A)$
is a completely isometric homomorphism onto a left ideal of
$\M_\ell(A)$.   Hence
$\lambda(A) = \M_\ell(A) \lambda(e)$.
Thus $\lambda(A)$ is a weak* closed subalgebra of $\M_\ell(A)$,
and so $B = \lambda(A)$ is a dual operator algebra.   If
we take a bounded net
$\lambda(a_i) \rightarrow \lambda(a)$ weak* in $\lambda(A)$,
then by definition of the weak* topology on $\M_\ell(A)$ from
\cite{BMDO} 3.2,
$a_i e = a_i  \rightarrow a e = a$ weak* in $A$.  Thus
$\lambda^{-1} $ is weak* continuous, so that by
Krein-Smulian (c.f. \cite{BMDO} Lemma 1.5)
$\lambda$ is weak* continuous.  
\end{proof}

\vspace{4 mm}
 
The following `1-sided version' of
Sakai's theorem may be known (certainly
most of it is contained in a result of Zettl \cite{Ze}
(recently given a new proof in
\cite{eor}, which one may greatly shorten in the
case where one has an operator space predual \cite{BMDO})).
 
\begin{theorem}  Let $J$ be a left ideal in a $C^*$-algebra,
and suppose that 
$J$  possesses a Banach space predual.  Then $M(J J^*)$ is a
$W^*$-algebra containing $J$ as a weak*-closed principal left ideal.
\end{theorem}
 
\begin{proof}  By the aforementioned  result of Zettl,
 $M(J J^*)$ is a $W^*$-algebra and $J$ is a dual operator space.
By
\ref{BLM1s}, $J$ has a
right identity $e$, so that $J = J J^* e = M(J J^*) e$.
\end{proof}

\vspace{4 mm}
         
Results such as \ref{upie} and \ref{Thlemi} are useful 
for deducing results about general  operator algebras with 
r.c.a.i., from results about left ideals in 
$C^*$-algebras.  Here is a sample application 
of this `transference
principle' (other 
examples will be given later):

\begin{corollary} \label{oares} 
 Let $A$ be an operator
 algebra with a right contractive
approximate identity, and also a right identity.   Then
$A$ has a right identity of norm 1, which is the limit 
in norm of the r.c.a.i.. 
\end{corollary}  

\begin{proof}  If $\{ e_\alpha \}$ is the r.c.a.i. for $A$,
then $\{ j(e_\alpha) \}$ is a r.c.a.i. for the left 
ideal envelope ${\mathfrak J}_e(A)$.
Similarly ${\mathfrak J}_e(A)$ and $A$ have a
common right identity.  Hence by \ref{babe}, our r.c.a.i.
converges in norm.
\end{proof}

\section{The Banach-Stone theorem}  

We prove several stages, or cases, of
this theorem, which asserts essentially that
linear surjective complete isometries between
left ideals of $C^*$-algebras (resp. between operator algebras
with r.c.a.i.),  are characterized as a
composition of a translation by a partial isometry $u$,
and a surjective completely isometric
homomorphism onto another right ideal (resp. operator algebra
with r.c.a.i.) which is
a translate of one of the original ideals (resp. algebras) by 
$u^*$. 
To see that the `translate by a partial isometry' is not artificial,
consider an infinite dimensional Hilbert space $H$ and $S$ the
shift operator.   Set $I = B(H)$ and $J = B(H) S$.  These
ideals are clearly linearly completely isometric, but
there is no homomorphism of $I$ onto $J$
(since $J$ has no 2-sided identity).   This example shows that
the following theorem (which comprises Case (1)) 
 is best possible:

\begin{theorem} \label{Th1}  Let $I$ and $J$ be principal left ideals
in $C^*$-algebras $A$ and $B$,
i.e.  $I = Ae$ and $J = Bf$, say, for orthogonal
projections $e, f$ in 
$I, J$ respectively.  Suppose also that 
$\varphi : I \rightarrow J$ is a linear surjective
complete isometry.  Then there exists a partial isometry
$u$ in $B$ with initial projection $f$, and a completely isometric 
surjective ideal homomorphism $\pi : I \rightarrow J_1$ 
such that $\varphi = \pi(\cdot) u$
and $\pi = \varphi(\cdot) u^*$.   Here $J_1 = B u^* = J u^* 
= B u u^*
\subset B$
is another left ideal of $B$ with right identity 
$u u^*$ of norm 1.  
 
Conversely, if $J$ is a left ideal of a 
$C^*$-algebra $B$, and if $u$ is a partial isometry in 
$B$ with initial projection a right identity for 
$J$ of norm $1$, then $J u^* = B u^* = B u u^*$ is a 
left ideal $J_1$ of $B$ with right identity
$u u^*$ of norm 1, and $J_1$ is 
linearly completely isometrically isomorphic to $J$ via
right multiplication by $u^*$.     Hence the composition
of right multiplication by $u^*$, with any 
completely isometric
surjective homomorphism $I \rightarrow J_1$, is a
linear completely isometric isomorphism $I \rightarrow J$.   
\end{theorem}

\begin{proof}   Recall from the introduction
that a completely 
isometric surjection between TRO's is a triple morphism.
Hence $\varphi$ is a 
triple isomorphism. 
 Therefore if $u = \varphi(e)$
then it is easy to check that $\pi(\cdot) =
\varphi(\cdot) u^*$ is a homomorphism
onto $J u^*$.   Similar considerations show that $p = u u^*$ is
an idempotent, which is an orthogonal projection since it is
selfadjoint.   Thus $u$ is a partial isometry.
We claim that $u^* u = f$.  To see this note 
that $u^* u$ is an orthogonal projection, and that for 
any $\varphi(x) \in J$ we have 
 $\varphi(x) u^* u = \varphi(x e) = \varphi(x)$, using the 
definition of a triple morphism.  Thus   
$f u^* u = f$.  On the other hand, $u^* u f = u^* u$ since 
$u \in Bf$.  Hence $f = u^* u$.   Also, 
$J u^* = B f u^* = B  u^* u u^* = B u^*$.   Defining $J_1$ to be 
this last
space we see that it is clearly a left ideal of $B$, and $J_1$ 
contains  $u u^*$, which is indeed a right identity of 
 norm 1 for $J_1$ since $u$ is a partial isometry.  Thus 
$J_1 = B u u^*$ too.

Since $\pi(\cdot) =
\varphi(\cdot) u^*$ we obtain 
$\pi(\cdot) u = \varphi(\cdot) u^* u = \varphi(\cdot) f 
= \varphi(\cdot)$.   It follows from this too that 
$\pi$ is a complete isometry, and therefore also a
triple morphism.  Thus $\pi$ is a completely isometric
ideal homomorphism.

Conversely, if $J, B, u$ are as stated, then 
$J = B u^* u$ so that $J u^* = B u^*$ which is also a 
right ideal of $B$.  Clearly the last space equals 
$ B u u^*$ since $B u u^* \subset B u^* = B u^* u u^*
\subset B u u^*$.  The remainder of 
the converse direction is left to the reader.  
\end{proof}

\vspace{4 mm}

Thus we see that central to the Banach-Stone theorem in
our 1-sided case, is a certain Murray - von Neumann equivalence
of projections.

Note in the `unital case' of the above theorem, i.e. the 
case that $\varphi(e) = f$, then $u = u^* = f$, $J_1 = J$, and
$\varphi$ itself is a homomorphism.  Conversely, if 
$\varphi$ is a homomorphism, then necessarily $\varphi(e) = f$,
as is easy to see from \ref{pid}.  

\vspace{4 mm}

Having thoroughly analyzed the Banach-Stone theorem in Case (1),
we now move to Case (2).  Here we look at linear 
completely isometric isomorphisms $\varphi : A \rightarrow B$
between operator algebras
with a right identity of norm 1.   In the assertions
in  the first paragraph of the statement of the theorem,
and in the proofs of these assertions, we regard
$A$ and $B$ as having been identified with 
subalgebras of ${\mathfrak J}_e(A)$ and ${\mathfrak J}_e(B)$ respectively
(see \ref{upie}).   Thus mention 
of the `canonical Shilov embedding homomorphisms'
$j$ have been suppressed, and all products and adjoints
in that paragraph
are taken in the containing $C^*$-algebra $\E(B) =
{\mathfrak J}_e(B) {\mathfrak J}_e(B)^*$.   

\begin{theorem} \label{Th2}  (Banach-Stone for operator algebras
with right identities.)
Suppose that $\varphi : A \rightarrow B$
is a surjective linear
completely isometric isomorphism between operator algebras
with a right identity of norm 1.  Then there exists a 
partial isometry $u \in {\mathfrak J}_e(B)$ (indeed in $B$)
with initial projection the 
right identity of $B$, such that the subspace $B' = B u^*$ 
of $\E(B)$ is a subalgebra (and consequently an
operator algebra) with a right identity $u u^*$ of norm 1;
and there exists a completely isometric surjective homomorphism 
$\pi : A \rightarrow B'$, such that $\varphi =
\pi(\cdot) u$ and $\pi = \varphi(\cdot) u^*$.   Also, $u^* B 
\subset B$.

Conversely, suppose we 
are given a partial isometry $u$ on a Hilbert 
space $H$, such that $u$ lies in a subalgebra $B \subset
B(H)$,
such that the initial projection of $u$
is a right identity of $B$, and such that 
$u^* B \subset B$.  Then $B' = B u^*$ is an
 operator algebra with right identity
$u u^*$ of norm $1$, and $B'$ is linearly completely isometrically
isomorphic to $B$ via right multiplication by $u$.  Thus the 
composition of right multiplication by $u^*$, with any
completely isometric
surjective homomorphism $A \rightarrow B'$, is a
linear completely isometric isomorphism $A \rightarrow B$.
\end{theorem}    

We remark that
in the language of \cite{BK}, $u$ and $u^*$ are in $LM(B)$.
  
\begin{proof}  
Suppose that $\varphi : A \rightarrow B$
is a linear
completely isometric isomorphism, and extend 
$\varphi$ to a linear
completely isometric isomorphism 
$\bar{\varphi} : {\mathfrak J}_e(A) \rightarrow 
{\mathfrak J}_e(B)$ (such extension exists 
by Hamana theory (\cite{Ham2} or \cite{BShi} Appendix 
A)).  By Theorem
 \ref{upie}, 
${\mathfrak J}_e(A)$ is a left ideal of the $C^*$-algebra $\E(A)$, and
${\mathfrak J}_e(A)$ has right identity $e$.  Similar assertions hold for 
${\mathfrak J}_e(B)$.  Thus by the proof of \ref{Th1}, if $u = \varphi(e) = 
\bar{\varphi}(e)$ then $u$ is a partial isometry in $B$, with 
$u^* \in B^* \subset \T(B)^* \subset \E(B)$, whose initial projection is
$f$, and $\pi = \bar{\varphi}(\cdot) u^*$ is a completely isometric
surjective homomorphism 
${\mathfrak J}_e(A) \rightarrow 
{\mathfrak J}_e(B) u^*$.  The restriction 
of $\pi$ to $A$ maps onto the subalgebra $B u^*$ of 
$\E(B)$.  Since $u$ is a partial isometry,
$u u^*$ is indeed a right identity of $B u^*$.   Finally, 
since $B u^* B u^* \subset B u^*$, post multiplying by $u$ 
gives $B u^* B \subset B$, so that
 $u^* B = u^* u u^* B = f u^* B \subset B u^* B \subset B$.

Conversely, given $u$ as stated, then since $u^* B \subset B$ we have
that $B u^*$ is a subalgebra of $\E(B)$ with right identity $u u^*$.
The remainder of the converse direction is obvious.
\end{proof}

\vspace{4 mm}
     
One can prove further that ${\mathfrak J}_e(B') = 
{\mathfrak J}_e(B) u^*$, and that 
$\E(B') = \E(B)$, but we omit the details.

\begin{corollary} \label{ca2co}  Suppose that $\varphi : A \rightarrow B$
is a surjective linear
completely isometric isomorphism between operator algebras with right 
identities $e$ and $f$ of norm $1$.  Then $\varphi$ is a homomorphism
if and only if $\varphi(e) = f$.
\end{corollary} \begin{proof}  The one direction follows from \ref{pid},
 the other by noting that if we follow the proof of \ref{Th2}, then
$\bar{\varphi}(e) = f$, so that $\bar{\varphi}$ is a homomorphism by
the remarks after \ref{Th1}.
\end{proof}

\begin{corollary} \label{ca2co3}  Suppose that $A$ is an 
operator algebra with a right identity of norm $1$, and 
suppose that $A$ has another product $m : A \times A \rightarrow A$  
with respect to which $A$ is completely isometrically isomorphic 
to an operator algebra with a right identity of norm $1$.
Then there is a partial isometry $u \in {\mathfrak J}_e(A)$ (and indeed in 
$A$) such that $m(x,y) = x u^* y$ for all $x,y \in A$. 
Indeed $u$ is the right identity for $m$, and $u^* u$ is
the right identity for the first product.
\end{corollary} 

\vspace{4 mm}
This result also has corollaries of the type given in 
\cite{Bcomm}; we leave the details to the reader. 
 
We now turn to Case (3) of the Banach-Stone theorem:

\begin{theorem} \label{bsgid}  (The Banach-Stone theorem
for left ideals in $C^*$-algebras.)  Consider
a surjective linear complete isometry $\varphi : I \rightarrow J$
between arbitrary
left ideals of $C^*$-algebras.   Let $\E = J J^*$, and let $\M$ 
be the von Neumann algebra $(J J^*)^{**}$.   
Then there exists another
left ideal $J_1$ of $\E$, with $J_1 J_1^* = \E$,
and a surjective completely isometric ideal homomorphism
$\pi : I \rightarrow J_1$.  Moreover there exists a partial isometry 
$u \in \M$ such that the initial projection of $u$ is 
the right identity of $J^{**}$ (indeed of $RM(J)$ - see \S 4),
and such that $J_1 = J u^*$, $J = J_1 u$, and such that
$\varphi = \pi(\cdot) u$, and $\pi = \varphi(\cdot)u^*$.
\end{theorem}

There clearly exists a (tidier) converse to this theorem, 
as in Cases (1) and (2).

\vspace{4 mm}

\begin{proof}
Consider the 
second dual  $\varphi^{**} : I^{**} \rightarrow J^{**} \subset \M$,
and now we are back in Case (1).   For if $I$  is a left ideal of a
$C^*$-algebra $A$, then $I^{**}$ is a left ideal of $A^{**}$, but 
now $I^{**}$ has a right identity $e$ of norm $1$, namely $e$
is a weak*-accumulation
point of the r.c.a.i. of $I$ (by \cite{BonsallDuncan} 28.7).   Thus by 
Case (1) we have that $u = \varphi^{**}(e)$ is a partial isometry in
$J^{**} \subset \M$, and the  
initial projection of $u$ is
the matching right identity of $J^{**}$.  Moreover
$\pi = \varphi^{**}(\cdot) u^*$ is a completely isometric homomorphism
and so on.  Restricting $\pi$ to $I$ 
gives a completely isometric homomorphism
$\pi'$ onto the subalgebra  $J_0 = J u^*$ of $\M$,
and $\varphi$ is the composition of $\pi'$
with a right translation by $u$.   Moreover, $\pi'$ is
easily seen to be a triple morphism:
$$\pi'(x) \pi'(y)^* \pi'(z) =
\varphi^{**}(\hat{x}) u^* u \varphi^{**}(\hat{y})^* \varphi^{**}(\hat{z})
u = \varphi(x) \varphi(y)^* \varphi(z) u 
=  \varphi(xy^*z) u = \pi'(xy^*z)$$
for $x,y,z \in I$.  Thus $\pi'$ is a completely isometric
ideal homomorphism.  Therefore by Lemma \ref{babi},
$\pi'$ is the restriction of a surjective 1-1
*-homomorphism  $I I^* \rightarrow
J_0 J_0^*$.  Thus $J_0 J_0^* = J u^* u J^* = \E$ contains 
$J_0$ as a left ideal; 
or to be more precise, $\hat{\E}$ contains $J_0$.
Thus we may regard $\pi'$ as a completely isometric homomorphism
$\pi : I \rightarrow 
J_1$ onto a right ideal $J_1$ of $\E$ (note
$\hat{J_1} = J_0$).   The rest is clear.
\end{proof}

\vspace{4 mm}
     
Finally, Case (4) of the Banach-Stone theorem, i.e. the case 
of a surjective linear complete isometry between arbitrary operator 
algebras with r.c.a.i..  Again it is clear that by passing to the 
second dual and using Case (2) in the way we tackled Case (3)
using Case (1), {\em or} using Case (3) in the way we 
tackled Case (2) using Case (1), will give a correct theorem
resembling Theorems \ref{Th1}, \ref{Th2}, and \ref{bsgid}.  We leave the
details to the reader.   This should also be linked to  
the multiplier algebras in \cite{BK}.  

\begin{corollary}
\label{unt}  Let $\varphi : A \rightarrow B$ be a
surjective linear complete isometry between left ideals
of $C^*$-algebras, or between operator
algebras with r.c.a.i..  Then $\varphi$ is a homomorphism 
if and only if there exists a 
r.c.a.i. $\{  e_\alpha \}$ for $A$ such that $\{ \varphi(e_\alpha) \}$
is a r.c.a.i. for $B$. 
\end{corollary}

\begin{proof}   If the latter condition holds then 
$\varphi^{**} : A^{**} \rightarrow B^{**}$ is a 
surjective linear complete isometry.  Let $E$ be a 
weak* limit point of $\{  e_\alpha \}$ in $A^{**}$, and since $\varphi^{**}$ is 
weak*-continuous, $\varphi^{**}(E)$ is a weak* limit point of 
$\{ \varphi(e_\alpha) \}$.  So we are in the situation of 
Corollary \ref{ca2co} by \ref{namo},
 so that $\varphi^{**}$ and consequently
$\varphi$ is a homomorphism.  The converse direction is easier.
\end{proof}

\section{$LM(A)$ for an algebra with left contractive approximate 
identity}

In this section we develop the 
`left multiplier operator algebra' $LM(A)$ of 
an operator algebra with l.c.a.i. (note that of course 
$RM(A)$ for an operator algebra with r.c.a.i. will have the 
almost identical, `other-handed version', of this
 theory).  On the other hand the left multiplier operator algebra
of an operator algebra with r.c.a.i. turns out to have a quite 
different theory, which is studied in the sequel \cite{BK}, and
which we will not mention again in the present paper.    

Actually, it appears at the
present time as if there 
may be two classes of `good candidates' for $LM(A)$
if $A$ has l.c.a.i..  
We will only really investigate the first of these classes in this
section, because in this case we can show that all the  candidates
in this class coincide, and we get a convincing theory
paralleling the known theory in many ways.  The one drawback
of this left multiplier algebra
 is that it does not contain the algebra itself
in general; but this is no surprise or mental obstacle
to anyone who has looked at the  `multiplier' or `centralizer'
 theory of nonunital Banach algebras.  
The second class of `candidates' consist of algebras which 
do contain $A$, but then of necessity one loses the 
important property that every nontrivial multiplier may be viewed as a
nontrivial map $A \rightarrow A$, and this introduces 
some problems.  Indeed this second class of 
`candidates' seems to be a rather unruly zoo; although 
there seems to be
a `best' or  `most canonical' candidate in this class, which
we call the `big left multiplier algebra' $BLM(A)$.
We will 
briefly discuss these matters
further at the end of this section (see 
also Remark 3 after \ref{remaf}).   

Thus we restrict our attention to the first class
of candidates.   Since this
follows closely the 
essentially  known theory for the case of a two-sided c.a.i.,
we will try to be brief.  As a 
historical note, the latter 
theory was begun in \cite{PuR}, and also independently 
developed around that time from a different angle, mostly by  
Paulsen \cite{Punp}, as a tool for our project 
on Morita equivalence 
\cite{BMP}, 
although this latter material was not circulated or
published.   Some other facts about this case of 
an operator algebra with a two-sided c.a.i. have arisen in 
several of the authors
 papers over the years, for example \cite{Bhmo,BShi}.
However since these facts were not particularly profound, 
followed mostly from the ideas of \cite{PuR} and \cite{Punp}, 
and did not play a crucial role, we did not give a complete 
development in these places, essentially leaving the needed 
details to the reader.  
The majority of these facts are contained 
in the following result, and its proof.
In view of this result, if $A$ is an operator algebra
with l.c.a.i. then we will write $LM(A)$ for 
any of the completely isometrically isomorphic
algebras of the theorem, and $\mu_A : A \rightarrow LM(A)$
for the canonical map (see also Remark 2 below).   
To understand this result better it is good to have a simple 
example in mind, such as 
$A = R_n$ (the subalgebra of $M_n$ supported on 
the first row).

\begin{theorem} \label{onsi}  Let $A$ be an operator algebra
with l.c.a.i..  Then the following operator algebras are all
completely isometrically isomorphic:
\begin{enumerate}
\item [(1)] $\{ x \in A^{**} : x \hat{A} \subset \hat{A} \}/ker 
\; q$ where $q $ is the canonical homomorphism into $CB(A)$,
\item [(2)] $\M_\ell(A)$,
\item [(3)] $CB_A(A) \; \; \; \; $ (right module maps).
\end{enumerate}                                                                        
Moreover $CB_A(A) = B_A(A)$ isometrically.  If $A$ has a two-sided
c.a.i., then $ker \; q = (0)$ in {\em (1)}, and also 
the above algebras 
also are completely isometrically isomorphic to:
\begin{enumerate} 
\item [(4)] the `maximum essential left multiplier extension
of $A$'.
 \end{enumerate}  
If $A$ satisfies condition
($\Li$) of \ref{defp} (for example if 
$A$ has a left identity of norm 1, or a two-sided c.a.i.,
or if $A$ is a right ideal of a $C^*$-algebra),
 then the algebras above are 
completely isometrically isomorphic to 
\begin{enumerate}
\item [(5)]   $\{ T \in B(H) : T \pi(A) \subset \pi(A) \}$, for
any completely $\Li$-isometric nondegenerate representation of $A$
(see definition after \ref{reps}),
\item [(6)]  $LM(B)$ where $B = \Li(A)$ (see \ref{defp}),
\item [(7)]  $\{ x \in B^{**} : x \hat{A} \subset \hat{A} \} \subset 
A^{**}$,
where $B = \Li(A)$.  
\end{enumerate}
\end{theorem}

To avoid distracting from other issues here,
we will not explain the term  `essential left multiplier extension' 
used in (4) until we reach it in the proof below.

\begin{proof}  
We first observe that for {\em any} operator algebra $A$
there are natural completely contractive
homomorphisms 
$$\{ x \in A^{**} : x \hat{A} \subset \hat{A} \} \rightarrow 
\M_\ell(A)  \rightarrow 
CB(A) .$$
Let us write $\sigma$ for the first homomorphism, and 
$\theta$ for the second.  From the `left handed variant' of
\ref{notfi}, the image of $\theta$ lies in $CB_A(A)$.
The fact that  $CB_A(A) = B_A(A)$ isometrically follows more
or less immediately from the relation $T(a) = \lim_\alpha T(e_\alpha) a$.
Next note that given $S \in CB_A(A)$,
then one may let $F$ be a weak* accumulation
point of $S(e_\alpha)$ in $A^{**}$, for the l.c.a.i.
$\{ e_\alpha \}$ for $A$.
Clearly $\Vert F \Vert
\leq \Vert S \Vert$.  For $a \in A$,
we have $$S(a) = \lim_\alpha S(e_\alpha a) =
\lim_\alpha S(e_\alpha) a = Fa \; . $$
Hence 
$q(F) = S$ where $q = \theta \circ \sigma$.  Thus $q$ is a 
quotient map, and similarly it is a complete quotient 
map.  Thus $\sigma$ is also a complete quotient map, and 
$ker \; \sigma = ker \; q$ since $\theta$ is 1-1.
This proves the completely isometric isomorphism between 
(1) and (3), and also between (1) and (2).   Thus $\M_\ell(A) \cong
CB_A(A)$ completely isometrically, which 
also shows that $CB_A(A)$ is a unital operator 
algebra
(or this fact may be proved directly).       

Now suppose that $A$ has property ($\Li$), and set $B = \Li(A)$
as in \ref{defp}.   Then $B^{**} \subset A^{**}$.
Examining the proof of (1) = (3) above, we see easily that the 
terms $S(e_\alpha)$ actually lie in $B$.  Hence the 
$F$ there lies in $\{ x \in B^{**} : x \hat{A} \subset \hat{A} \}$.
Thus the map $q$ mentioned above, restricted to the last set, 
is a complete quotient map too.  Therefore it is a 
complete isometry if we can show that it
 is 1-1.  To see this suppose 
that $F$ is in the set in (7) and $q(F) = 0$.  Then 
$F \hat{e_\alpha} = 0$.   This implies that $F = 0$, using the 
fact that a weak* limit point of the $\hat{e_\alpha}$ is a 
2-sided identity for $B^{**}$, and the fact
that the multiplication 
in a dual operator algebra is separately weak* continuous.
Thus we have that (3) = (7)  completely isometrically.
Note too that if $A$ is an operator algebra with 2-sided
c.a.i. then this shows that $ker \; q = (0)$ in (1).
Note that if $F$ is in the set in (7), then $F B \subset B$
quite clearly.  Conversely if $F B \subset B$ then  for $a \in A$
we have $F a = \lim F e_\alpha a \in A$ 
since $F e_\alpha a \in B a \subset A$.
This shows that (6) = (7).

Returning to (4), we need a definition of an
 `essential left multiplier extension' of 
an operator algebra $A$ with c.a.i..
For the purposes of this theorem we will define this to be 
a pair $(B,\pi)$ consisting of an operator algebra $B$ with
identity of norm 1, together with a
completely isometric homomorphism $\pi :
A \rightarrow B$, such that $\pi(A)$ is a left ideal of 
$B$, and such that the natural map $B \rightarrow CB(A)$ is 
completely isometric.   There is a natural ordering, and 
corresponding notion 
of equivalence, of `left multiplier extensions' of $A$, which we
will not bother to spell out.  However it is clear that since 
the just mentioned `natural map' maps into $CB_A(A)$, the 
algebras in (1)-(3) are the biggest essential left multiplier 
extensions. 
 
Finally to prove that  (5) = (6), we may w.l.o.g., 
by the definition
after \ref{reps} and the last assertion of
\ref{reps}, assume that $B = A$ is an operator algebra 
with 2-sided c.a.i., and that $\pi : A \rightarrow B(H)$ is 
a nondegenerate completely isometric homomorphism.  This 
case is well known; briefly, one way to see it is
as follows.  
Noting that the algebra in (5) is then an 
essential left multiplier extension of $A$, we see that the 
algebra in (5) is completely isometrically contained as a 
subalgebra of $CB_A(A)$.  Conversely, if $R \in B_A(A)$, we 
obtain a related map $T \in B(H)$ which may be defined by
$T \pi(a) \zeta = \pi(Ta) \zeta$, for $a \in A, \zeta \in H$.
Another way to see this quickly is using the well known fact 
that in this case, $H \cong A \hat{\otimes}_A H$.   We omit
the simple details, which as we said at the beginning of this
section, are essentially well known to experts.
\end{proof} 

\begin{corollary} \label{remaf} 
Let $A$ be an operator algebra with left identity
$e$ of norm 1.  Then $LM(A) = Ae$, which is a
unital subalgebra of $A$.  It is also a 
unital subalgebra of $\E(A)$, 
and $\E(A)$ is a unital $C^*$-algebra.
\end{corollary}

\begin{proof}  In this case $\Li(A) = Ae$ which is a
unital algebra.  Thus the first result here
follows 
from (6) of the previous theorem.  We saw in \ref{upie}
that $J = {\mathfrak J}_e(A)$ 
is a right ideal of a 
$C^*$-algebra, and 
$J$ has a left identity $e$.  Thus $\E(A) = J J^*$
has $e$ as a 2-sided identity.  Finally $A e \subset J J^* = \E(A)$.
\end{proof}

\vspace{4 mm}

\noindent {\bf Remark 1.}  There is a condition similar to (4) 
in \ref{onsi} which
is equivalent to (1)-(3) in general.  However it was sufficiently
more complicated to offset the benefits
of mentioning it in more detail.
                            
\vspace{4 mm}
               
\noindent {\bf Remark 2.}   It is very important that any kind of 
multiplier algebra $D$ of an algebra
$A$ should not merely be regarded as an algebra, 
but rather as a {\em pair} $(D,\mu)$ consisting also of a
homomorphism  
$\mu : A \rightarrow D$.   Sometimes we write 
$\mu_A$ to indicate the dependence on $A$.  Saying that 
two such algebras $(D,\mu)$ and $(D',\mu')$ 
are the same as multiplier algebras 
must entail a completely isometric surjective homomorphism
$\theta : D \rightarrow D'$ such that 
$\theta \circ \mu = \mu'$.  In this case we say that 
$D$ and $D'$ are $A$-isomorphic.  
Thus in each of the seven equivalent 
formulations of the previous theorem, we need to have in mind 
also what the map $\mu$ is in each case.  In (1) it is the map $a \mapsto
\hat{a} + \; ker \; q$; in (2) and (3) it is essentially the 
regular representation $\lambda$; in (5) it is $\pi$; in 
(6) it is the natural left representation of $A$ on its 
left ideal $\Li(A)$; and in (7) it is $a \mapsto \hat{a} E$,
where $E$ is as in the remark before \ref{reps}.
All these maps 
are completely contractive homomorphisms.  One then needs 
to check that these seven algebras are all $A$-isomorphic.
We leave these assertions 
to the reader who wishes to be more careful.

\vspace{4 mm}
     
\noindent {\bf Remark 3.}   Suppose that $A$ is an operator algebra
with l.c.a.i., and that $\pi : A \rightarrow B(H)$ is a completely 
isometric representation.   Define
$LM(\pi) = \{ T \in B(H) : T \pi(A) \subset \pi(A) \}$,
the left idealizer of $\pi(A)$ in $B(H)$.
Then it is straightforward to exhibit a 
completely contractive homomorphism $\sigma : 
LM(\pi) \rightarrow LM(A) = CB_A(A)$.  Conversely,
given $T \in  CB_A(A)$, taking a weak operator 
limit point $S$  of $\pi(T(e_\alpha))$ gives 
$S \in LM(\pi)$.  This is really saying that 
$LM(A) \cong \; LM(\pi)/ker \; \sigma$ completely 
isometrically isomorphically.  One may view this 
observation as an attempt to remove the use of property
($\Li$) in (5).   
Perhaps an
investigatiion of this quotient might be tied
to Sarason's semi-invariant subspace technique
(see \cite{Arv1} for example).  It is interesting to 
note that if $\pi$ is the usual representation of 
$R_2$, then $LM(\pi)$ is a 3-dimensional operator 
algebra (this was pointed out to me by M. Kaneda).    
 Note that $LM(\pi)$ is
highly dependent on $\pi$, to see this consider 
$R_2$ again; the natural 
representation $\pi$ has $LM(\pi)$ 3-dimensional.
However, if $\sigma = \pi \oplus \epsilon$, where
$\epsilon$ is the projection onto the 1-1 coordinate,
then $LM(\sigma)$ is strictly larger.  It would be interesting
 to see if there is a nonrestrictive condition
under which one obtains `independence from the 
particular $\pi$ used'.
 
\vspace{4 mm}
           
We now turn to the notion which in the $C^*$-algebra literature 
is referred to as `essential homomorphisms' or sometimes 
`nondegenerate homomorphisms'.   For our purpose we shall use the
name `$A$-nondenerate morphism'.   For us this shall mean 
a completely contractive homomorphism $\pi : A \rightarrow LM(B)$ 
satisfying the following equivalent conditions:

\begin{theorem} \label{Anon}  Let $A$ and $B$ be two operator 
algebras with l.c.a.i.'s, and let $\pi : A \rightarrow LM(B)$ 
be a completely contractive homomorphism.    The following 
are equivalent:
\begin{itemize}
\item [(i)]  There exists a l.c.a.i. $\{ e_\alpha \}$
for $A$ such that $\pi(e_\alpha) b \rightarrow b$ for all 
$b \in B$,
\item [(ii)]  For every l.c.a.i. $\{ e_\alpha \}$ for $A$,
we have $\pi(e_\alpha) b \rightarrow b$,
\item [(iii)]  $B$ is a nondegenerate left $A$-module via $\pi$,
\item [(iv)]  Any $b \in B$ may be written $b = \pi(a) b'$ 
for some $a \in A, b' \in B$.
\end{itemize}
If these conditions hold, there exists a completely contractive 
unital homomorphism
$\hat{\pi} : LM(A) \rightarrow LM(B)$  such that 
$\hat{\pi} \circ \mu_A = \mu_B$.  Moreover 
$\hat{\pi}$ has the
property that $\hat{\pi}(x)(\pi(a) b) = \pi(xa) b$ for 
$x \in LM(A), a \in A, b \in B$;
 and it is the unique such homomorphism with this 
property.  Finally $\hat{\pi}$ is completely isometric 
if $\pi$ is completely isometric.
\end{theorem}

\begin{proof}   Clearly (i) implies that 
the span of terms $\pi(a) b$ is dense in $B$, which is
what we mean by nondegenerate.  So (i) implies (iii).  Clearly 
(iii) implies (ii), and (ii) implies (i), and (iv) implies (iii).
That (iii) implies (iv) follows from \cite{Pal} \S 5.2.

If these conditions hold, view $LM(A)$ and $LM(B)$ as in 
(3) of \ref{onsi}, and note the formula for 
$\mu_A$ in the previous remark.    We may follow the proof 
of Theorem 6.2 
in \cite{Bhmo} (which the author proved inspired 
by an argument of van Daele mentioned there).
The main difference is that we ignore the 
element $e$ mentioned there, which we can get away with 
by taking $d$ there to be the l.c.a.i. from $A$.
One also needs to use \cite{Pal} 5.2.2, and the matricial 
version of it, in order to show that $\hat{\pi} : LM(A)
\rightarrow CB_B(B)$, and that it is a complete 
contraction.   The remaining problem is the last assertion.
Supposing that $\pi$ is  completely isometric, 
we have $$\Vert \hat{\pi}(T) \Vert_{cb} \geq
\Vert [\hat{\pi}(T) (\pi(a_{ij}) b_{kl})] \Vert  
= \Vert [\pi(T a_{ij}) b_{kl}] \Vert$$
providing that $\Vert [a_{ij}] \Vert, \Vert [b_{kl}] \Vert
\leq 1$.  Taking the supremum over all such 
$[b_{kl}] \in M_m(B)$, gives that 
$$\Vert \hat{\pi}(T) \Vert_{cb} \geq    \Vert [\pi(T a_{ij})] \Vert
= \Vert [T a_{ij}] \Vert \; \; .  $$  
Taking the supremum over all such $[a_{ij}] \in M_n(A)$
gives that $\Vert \hat{\pi}(T) \Vert_{cb} \geq
\Vert T \Vert_{cb}$.  So $\hat{\pi}$ is isometric,
and similarly it is completely isometric.     
\end{proof}   

\vspace{4 mm}

We note that the canonical map $\mu_A : A \rightarrow LM(A)$ 
is an $A$-nondegenerate morphism.  Note too that 
if $\pi : A \rightarrow B$ is a completely contractive 
homomorphism between  two operator
algebras with l.c.a.i.'s, then we can say that 
$\pi$ is $A$-nondegenerate, if $\mu_B \circ \pi$ 
is $A$-nondegenerate.  This happens if and only if 
$\pi$ takes a l.c.a.i. for $A$ to a l.c.a.i. for $B$, 
or equivalently that $B$ is a nondegenerate $A$-module 
via $\pi$.   Thus the inclusion map 
from a closed  subalgebra $A$ of an 
operator algebra $B$ is $A$-nondegenerate
if  $A$ contains a l.c.a.i. for $B$.  In this case
we have

\begin{corollary}  \label{incmu}  Let $A$ be a 
closed  subalgebra of an
operator algebra $B$, and suppose that $A$ contains a 
l.c.a.i. for $B$.  Then $LM(A) \hookrightarrow LM(B)$ 
completely isometrically as a subalgebra.
\end{corollary}
      
\vspace{4 mm}

We end with the promised remarks about a second class
of candidates for the `left
multiplier operator algebra' of an 
operator algebra $A$ with l.c.a.i..  
For a candidate $(B,\nu)$ in this class one would like to
have the property that $\nu : A \rightarrow B$
is a completely isometric homomorphism, such that 
$B \nu(A) \subset \nu(A)$.  
Unfortunately  then 
one must lose the useful `essential'
condition (i.e. that $x \nu(A) = 0$ implies $x = 0$).
Thus for many $A$ (such as our standard example 
$R_n$) a subset of $B$ plays no role in any action on
$A$, and this fact does not seem to
bode well for  `uniqueness'
properties of such multiplier algebras.   Indeed this takes us
out of the classical `multiplier'/`centralizer' framework 
from Banach algebra theory (\cite{Pal} \S 1.2 for example), where
a multiplier which annihilates $A$ must be the zero multiplier.   
Also one cannot hope for conditions like (1)-(3) of 
\ref{onsi}, it appears.  

With these cautions it seems nonetheless that there  is a
`most interesting  candidate' in this class.
Consider the injective $C^*$-algebra
 $B$ and homomorphism $j$ of (the other-handed version of) 
Theorem \ref{Thlemi},
and define 
 $BLM(A) = \{ x \in B :
x j(A) \subset j(A) \}$.  That is, $BLM(A)$ is the left idealizer 
of $j(A)$ in $B$.   
For example, $BLM(R_2)$ is the
upper triangular $2 \times 2$ matrix algebra. 
Since $I_{11} = eBe \subset
B = I_{22}$ in the language of
(the other-handed version of) \ref{Thlemi}, 
it may be proved by the 
multiplication theorem in \cite{BPnew} and 1.3 there, that 
$BLM(A)$ is the `biggest' 
unital operator algebra containing 
$A$ completely isometrically as a subalgebra, in the sense 
that it contains a completely contractive image of every 
other such algebra $C$.  It may not contain $C$ itself though
(consider $A = R_2$, and see Remark 3 above).    
    
Other candidate definitions for a left multiplier algebra containing $A$
might look like the $LM(\pi)$ algebras in Remark 3 
after \ref{remaf}.  One might hope that there is a suitable and
not too strong condition on $\pi$ there so that these
algebras are independent of the particular $\pi$, but 
this seems unlikely at present.   
It seems as if one could present large numbers of further 
such candidates, which at this point 
seem unrelated in any way.


\begin{thebibliography}{99}
 
\bibitem{Arv1}  W. B. Arveson,
{\em Subalgebras of }$C^{*}-${\em algebras,}
 Acta Math. {\bf 123 }(1969), 141-224.

\bibitem{Arv2}   W. B. Arveson,
{\em Subalgebras of }$C^{*}-${\em algebras II,}
Acta Math. {\bf 128} (1972), 271-308.

\bibitem{Bcomm}  D. P. Blecher,
{\em  Commutativity
in operator algebras,} {\em Proc. Amer.  Math. Soc.}
{\bf 109} (1990), 709-715.
                      
\bibitem{Bhmo} D. P. Blecher,
 {\em  A generalization of Hilbert modules,}
J. Funct. An. {\bf 136} (1996), 365-421.    
                         
\bibitem{Bnat}  D. P. Blecher,  {\em Some general theory of
operator algebras and their modules,} in
 {\em Operator algebras and
applications\/},
A. Katavalos (editor), NATO ASIC, Vol.
495, Kluwer, Dordrecht, 1997.

\bibitem{BShi}  D. P. Blecher,  {\em  The Shilov boundary
of an operator space, and the characterization
theorems,}  J. Funct. An. {\bf 182} (2001), 280-343.

\bibitem{BMDO}  D. P. Blecher,  {\em Multipliers and dual operator algebras,}
J. Funct. An.   {\bf 183} (2001), 498-525.
 
\bibitem{BEZ}  D. P. Blecher, E. G.
Effros and V. Zarikian, {\em One-sided M-ideals and multipliers 
in operator spaces I,} Preprint (2000).  To appear {\em  Pacific J. Math}.     

\bibitem{BK}  D. P. Blecher and M. Kaneda, {\em 
The ideal envelope of an
operator algebra}, Preprint.

\bibitem{BMP} D. P. Blecher, P. S. Muhly and V. I. Paulsen,  
{\em  Categories of
operator modules - Morita equivalence and projective modules,}
 {\em Memoirs of the American Mathematical Society}
Vol. 143, number 681 (Jan, 2000). 

\bibitem{BPnew}  D. P. Blecher and V. I. Paulsen,
 {\em  Multipliers
of operator spaces, and the injective envelope,}
 Pacific J. Math.    {\bf 200} (2001), 1-17.  

\bibitem{BRS}  D. P.  Blecher, Z. J. Ruan, and A. M.
Sinclair, {\em A characterization
of operator algebras}, J. Funct. An. {\bf 89} (1990),
188-201.

\bibitem{BonsallDuncan}  F. F. Bonsall and J. Duncan,
{\em Complete normed
algebras,}      Springer-Verlag, New York-Heidelberg
(1973).                                                                                         
\bibitem{Con}  J. B. Conway, {\em A course in functional
analysis,} Graduate Texts in Mathematics, Springer-Verlag,
1985.

 
\bibitem{eor}   E. G.
Effros, N. Ozawa and Z. J. Ruan, {\em On injectivity and
nuclearity for operator spaces,} Duke Math. J. To appear.

\bibitem{ERns}  E. G. Effros and Z. J. Ruan,
{\em On non-self-adjoint operator
algebras},  Proc. Amer. Math. Soc. {\bf 110} (1990),
915-922.     
     
\bibitem{ERbook}   E. G.
Effros and Z. J. Ruan, {\em Operator Spaces}, Oxford University
Press (2000).                                                             

\bibitem{ECP} K. El Amin, A. M. Campoy, and
A. Rodriguez-Palacios, {\em
A holomorphic characterization of $C^*$- and JB$^*$-algebras,}
Manuscripta Math. {\bf 104} (2001), 467-478.


\bibitem{Fill}  P. A. Fillmore, {\em A users guide to 
operator algebras,} Canadian Math. Soc. series,
Wiley-Interscience (1996).    

 \bibitem{Ham1}  M. Hamana, {\em    Injective
envelopes of operator systems,}
Publ. R.I.M.S. Kyoto Univ. {\bf 15} (1979), 773-785.
 
\bibitem{Ham2}  M. Hamana, {\em Triple  envelopes and Silov
boundaries of operator spaces,}  Math. J. Toyama University
{\bf 22} (1999), 77-93.   

\bibitem{Ka}  R. V. Kadison, {\em Isometries of
operator algebras,}  Ann. of Math. {\bf 54} (1951),
325-338.    

\bibitem{KP}  M. Kaneda and  V. I. Paulsen, {\em
Characterizations of essential ideals as operator modules
over $C^*$-algebras}, Preprint (2000).


\bibitem{LM3}  C. Le Merdy, {\em An operator space
characterization of dual operator algebras,}  Amer. J. Math.,
{\bf 121} (1999), 55-63.
               
 
\bibitem{NR}  M. Neal and B. Russo, {\em A holomorphic characterization
of ternary rings of operators}, Preprint, 2001.

\bibitem{Pal}  T. W. Palmer, {\em Banach Algebras and the general 
theory of *-algebras Vol. I,} Encyclopedia of mathematics and its
applications Vol. 49, Cambridge University Press (1994).

\bibitem{P}   V. I. Paulsen, {\em Completely bounded maps 
and dilations,} Pitman Research
Notes in Math., Longman, London, 1986.

\bibitem{Punp}  V. I. Paulsen, {\em Unpublished notes on 
multiplier algebras,} (circa 1992).

\bibitem{Ped} G. Pederson, {\em C$^*$-algebras
and their automorphism groups,}
Academic Press (1979).                         
 
\bibitem{Pis} G. Pisier, {\em Introduction to the theory
of operator spaces,} to appear.

\bibitem{PuR}  Y-t. Poon and Z. J. Ruan, {\em Operator algebras with
contractive approximate identities, } Canadian J. Math. {\bf 46}(1994),
397-414.     

\bibitem{W-O}  N. E. Wegge-Olsen, {\em K-theory and $C^*$-algebras,}
 Oxford Univ. Press (1993).
                          
\bibitem{Ze}  H. Zettl, A characterization of ternary rings of operators,
{\em Adv. Math.} {\bf 48} (1983), 117--143.
                         
\end{thebibliography}
\end{document}